\documentclass[12pt]{article}
\usepackage{graphicx}
\usepackage{amsmath,amssymb,amsthm,amsfonts}
\newtheorem{thm}{Theorem}[section]
\newtheorem{cor}[thm]{Corollary}

\newcommand{\R}{{\mathbb{R}}}

\newcommand{\1}{\partial}

\begin{document}
\title{Gradient estimates for a nonlinear parabolic equation under Ricci flow}
\author{Shu-Yu Hsu\\
Department of Mathematics\\
National Chung Cheng University\\
168 University Road, Min-Hsiung\\
Chia-Yi 621, Taiwan, R.O.C.\\
e-mail:syhsu@math.ccu.edu.tw}
\date{June 25, 2008}
\smallbreak \maketitle
\begin{abstract}
Let $(M,g(t))$, $0\le t\le T$, be a n-dimensional complete noncompact 
manifold, $n\ge 2$, with bounded curvatures and metric $g(t)$ evolving by 
the Ricci flow $\frac{\partial g_{ij}}{\partial t}=-2R_{ij}$. We will extend 
the result of L.~Ma and Y.~Yang and prove a local gradient estimate for 
positive solutions of the nonlinear parabolic equation 
$\frac{\1 u}{\1 t}=\Delta u-au\log u-qu$ where $a\in\R$ is a constant
and $q$ is a smooth function on $M\times [0,T]$. 
\end{abstract}

\vskip 0.2truein



\vskip 0.2truein
\setcounter{equation}{0}
\setcounter{section}{0}

\section{Introduction}
\setcounter{equation}{0}
\setcounter{thm}{0}

Recently there is a lot of studies on Ricci flow on manifolds by R.~Hamilton
and others \cite{H}, \cite{CLN}, because Ricci flow is a powerful tool in 
analyzing the structure of manifolds. Let $(M,g(t))$, $0\le t\le T$, be 
a n-dimensional complete noncompact manifold with metric $g(t)$ evolving 
by the Ricci flow 
\begin{equation}
\frac{\1 g_{ij}}{\1 t}=-2R_{ij}\quad\mbox{ in }M.
\end{equation}
$(M,g(t))$ is said to be a gradient Ricci soliton (\cite{H},\cite{T}) 
if there exists a constant $c\in\R$ such that
\begin{equation}
\mbox{Ric}=cg+D^2f\quad\mbox{ in }M
\end{equation}
where $D^2f$ is the Hessian of $f$. $(M,g(t))$ is a expanding gradient 
soliton if $c<0$ and $(M,g(t))$ is a shrinking gradient soliton if $c>0$.
As observed by L.~Ma \cite{M} if we let $f=\log u$, then after some 
computation using (1.2) we get
\begin{equation}
\Delta u+2cu\log u=(A_0-nc)u\quad\mbox{ in }M
\end{equation}
for some constant $A_0$. Then the study of the gradient soliton on $M$ is 
reduced to the study of the properties of the solution $u$ of (1.3).
In \cite{M} L.~Ma using the technique of S.Y.~Cheng, P.~Li and S.T.~Yau
\cite{CY}, \cite{LY}, proved the local gradient estimates of positive
solutions of the equation
\begin{equation*}
\Delta u-au\log u-bu=0\quad\mbox{ in }M
\end{equation*}
where $a>0$ and $b\in\R$ are constants for complete noncompact manifolds 
with a fixed metric and curvature locally bounded below. In \cite{Y} 
Y.~Yang extended L.~Ma's result and proved the local gradient estimates of
positive solutions of the equation
\begin{equation*}
u_t=\Delta u-au\log u-bu\quad\mbox{ in }M\times (0,T]
\end{equation*}
where $a,b\in\R$ are constants for complete noncompact manifolds with a fixed
metric and curvature locally bounded below. In \cite{CTY}
A.~Chau, L.F.~Tam and C.~Yu, proved the local gradient estimates of 
positive solutions of the conjugate heat equation
\begin{equation*}
u_t=\Delta^tu-qu\quad\mbox{ in }M\times (0,T]
\end{equation*}
where $q$  is a smooth function on $M\times [0,T]$ and $\Delta^t$ is the 
Laplace operator with respect to the metric $g(t)$ for any complete 
noncompact manifold $(M, g(t))$, $0\le t\le T$, with bounded curvature 
and metric $g$ satisfying
\begin{equation}
\frac{\1 g_{ij}}{\1 t}=2h_{ij}\quad\mbox{ in }M\times (0,T]
\end{equation}
where $h(x,t)=(h_{ij}(x,t))$ is a smooth symmetric tensor on $M\times [0,T]$. 
In this paper we will extend their results and prove the local gradient 
estimates of positive solutions of the equation
\begin{equation}
u_t=\Delta^tu-au\log u-qu\quad\mbox{ in }M\times (0,T]
\end{equation}
where $a\in\R$ is a constant and $q$ is a smooth functions on $M\times [0,T]$ 
for any complete 
noncompact manifold $(M, g(t))$, $0\le t\le T$, with bounded curvature 
and metric $g$ satisfying (1.4).  In particular
our result holds for complete noncompact manifolds with metric evolving by 
Ricci flow.

We start with some definitions. Let $\nabla^t$ and $\Delta^t$ be the 
covariant derivative and Laplacian with respect to the metric $g(t)$. 
When there is no ambiguity, we will drop the superscript and write $\nabla$,
$\Delta$, for $\nabla^t$, $\Delta^t$, respectively. For any $r>0$, 
$x_0\in M$, $0<t_1\le t_0\le T$, let $B_r(x_0)$ be the 
geodesic ball with center $x_0$ and radius $r$ with respect to the 
metric $g(0)$.
For any $x_1,x_2\in M$, let $r(x_1,x_2)$ be the distance between $x_1$ 
and $x_2$ with respect to the metric $g(0)$.

\section{Local gradient estimates}
\setcounter{equation}{0}

In this section we will prove the local gradient estimates for positive
solutions of (1.5). We will let $(M, g(t))$, $0\le t\le T$, be a complete 
noncompact manifold with metric $g$ satisfying (1.4) and 
$$
|Rm|\le k_0\quad\text{ on }M\times [0,T]
$$
for some constant $k_0>0$ for the rest of the paper. Similar to \cite{CTY} 
we will assume that 
$q$ is a smooth functions on $M\times [0,T]$ and $h=(h_{ij})$ is a 
smooth symmetric tensor on $M\times [0,T]$. We will also assume the following 
for the metric $g$ and symmetric tensor $h(x,t)=(h_{ij}(x,t))$:

\begin{enumerate}\rm
\item[(A1)] $\|h\|, \|\nabla^th\|$ are uniformly bounded on $[0,T]$ where
the norms are taken on $M$.  
\item[(A2)] $|q|,|\nabla^tq|, |\Delta q|$ are uniformly bounded
on $M\times [0,T]$.
\end{enumerate} 

We also let $H=g^{ij}h_{ij}$ be the trace of $h=(h_{ij})$.

\begin{thm}
Let $u$ be a positive solution of (1.5). Then for any $\alpha>1$ and 
$0<\delta<1$ 
there exists a constant $C_1>0$ depending on $k_0$, $\alpha$, $\delta$, 
$a_+=\max (a,0)$ and the space-time uniform bound of $\|h\|, \|\nabla^th\|$, 
$|q|$, $|\nabla^tq|$, $|\Delta^tq|$ such that
\begin{equation}
\frac{|\nabla u|^2}{u^2}-\alpha\frac{u_t}{u}-\alpha a\log u-\alpha q
\le\frac{n\alpha^2}{2(1-\delta)t}(1+ta_-)+C_1\quad\text{ in }M\times (0,T]
\end{equation}
where $a_-=-\min (a,0)$.
\end{thm}
\begin{proof}
We will use a modification of the technique of \cite{CTY} and \cite{Hs2} to 
prove the theorem. Let $f=\log u$. Then by (1.5) $f$ satisfies
\begin{equation}
f_t=\Delta f+|\nabla f|^2-af-q.
\end{equation} 
By \cite{LY} and \cite{CTY} $f$ satisfies
\begin{equation}
\left\{\begin{aligned}
(\Delta f)_t=&\Delta f_t-2h_{ij}f_{ij}-2h_{ik;i}f_k+\nabla H\cdot\nabla f\\
(|\nabla f|^2)_t=&2\nabla f_t\cdot\nabla f-2h(\nabla f,\nabla f)
\end{aligned}\right.
\end{equation}
where $h(\nabla f,\nabla f)=h_{ij}f_if_j$. Let
\begin{equation}
F(x,t)=t[|\nabla f|^2-\alpha (f_t+af+q)].
\end{equation}
By (2.3) and the Weitzenbock-Bochner formula \cite{A},
$$
\Delta (|\nabla f|^2)=2\nabla f\cdot\nabla (\Delta f)+2|\nabla^2f|^2
+2R_{ij}f_if_j,
$$
we have in normal coordinates
\begin{align}
t^{-1}\Delta F=&\Delta(|\nabla f|^2)-\alpha\Delta f_t
-\alpha a\Delta f-\alpha\Delta q\nonumber\\
=&2|\nabla^2f|^2+2\nabla f\cdot\nabla (\Delta f)+2R_{ij}f_if_j
-\alpha(\Delta f)_t-2\alpha h_{ij}f_{ij}-2\alpha h_{ik;i}f_k
\nonumber\\
&\qquad +\alpha\nabla H
\cdot\nabla f-\alpha a\Delta f-\alpha\Delta q.
\end{align}
Now by (2.2) and (2.4),
\begin{equation}
\Delta f=f_t+af+q-|\nabla f|^2
=\frac{1}{\alpha}\left(|\nabla f|^2-\frac{F}{t}\right)-|\nabla f|^2
=-\frac{1}{\alpha}\left((\alpha -1)|\nabla f|^2+\frac{F}{t}\right).
\end{equation}
Then by (1.4), (2.2), (2.4) and (2.6),
\begin{align}
&-\alpha(\Delta f)_t+2\nabla f\cdot\nabla(\Delta f)\nonumber\\
=&\frac{\1}{\1 t}\left(\frac{F}{t}+(\alpha -1)|\nabla f|^2\right)
+2\nabla f\cdot\nabla(f_t+af+q-|\nabla f|^2)\nonumber\\
=&\frac{F_t}{t}-\frac{F}{t^2}+2\alpha\nabla f\cdot\nabla f_t
+2\nabla f\cdot\nabla(af+q-|\nabla f|^2)-2(\alpha-1)h_{ij}f_if_j
\nonumber\\
=&\frac{F_t}{t}-\frac{F}{t^2}+2\alpha\nabla f\cdot\nabla\biggl (
\frac{1}{\alpha}|\nabla f|^2-\frac{F}{\alpha t}-af-q\biggr )
+2\nabla f\cdot\nabla(af+q-|\nabla f|^2)\nonumber\\
&\qquad -2(\alpha-1)h_{ij}f_if_j\nonumber\\
=&\frac{F_t}{t}-\frac{F}{t^2}-\frac{2}{t}\nabla f\cdot\nabla F
-2(\alpha -1)a|\nabla f|^2-2(\alpha -1)\nabla f\cdot\nabla q
-2(\alpha -1)h_{ij}f_if_j.
\end{align}
Hence by (2.5), (2.6) and (2.7),
\begin{align}
&t^{-1}(\Delta F-F_t+2\nabla f\cdot\nabla F)\nonumber\\
=&2\sum_{i,j}(f_{ij}^2-\alpha h_{ij}f_{ij})
+2R_{ij}f_if_j-(\alpha -1)a|\nabla f|^2-2(\alpha -1)\nabla f\cdot\nabla q
\nonumber\\
&\qquad -2(\alpha -1)h_{ij}f_if_j-\frac{F}{t^2}+a\frac{F}{t}
-2\alpha h_{ik;i}f_k+\alpha\nabla H\cdot\nabla f-\alpha\Delta q.
\end{align}
By (2.8) and Young's inequality, for any $0<\delta<1$
there exist constants $C_1>0$, $C_2>0$, $C_3>0$ and $C_4>0$ 
which depend on $a_+$ and such that
\begin{align}
&t^{-1}(\Delta F-F_t+2\nabla f\cdot\nabla F)\nonumber\\
\ge&2(1-\delta)\sum_{i,j}f_{ij}^2-C_1(|\nabla f|+|\nabla f|^2)
-C_2-\frac{F}{t^2}+a\frac{F}{t}\nonumber\\
\ge&\frac{2(1-\delta)}{n}(\Delta f)^2-C_3|\nabla f|^2
-C_4-\frac{F}{t^2}+a\frac{F}{t}.
\end{align}
By (1.4) and (A1) there exist constants $c_1>0$, 
$c_2>0$, such that
\begin{equation}
\left\{\begin{aligned}
&c_1g_{ij}(x,0)\le g_{ij}(x,t)\le c_2g_{ij}(x,0)\quad\forall 
x\in M, 0\le t\le T\\
&c_1g^{ij}(x,0)\le g^{ij}(x,t)\le c_2g^{ij}(x,0)\quad\forall 
x\in M, 0\le t\le T.\end{aligned}\right.
\end{equation}
Let $x_0\in M$ and $r(x)=r(x,x_0)$. Then by (2.10) there exists a constant
$c_3>0$ such that 
\begin{equation}
|\nabla^tr(x)|^2\le c_3\quad\mbox{ in }M\quad\forall 0\le t\le T.
\end{equation}
Let $\psi\in C^{\infty}([0,\infty))$, $0\le\psi\le 1$, be such that 
$\psi (r)=1$ for all $0\le r\le 1$, $\psi (r)=0$ for all $r\ge 2$ and 
$\psi''\ge -c_4$, $0\ge\psi'\ge -c_4$ and $\psi'{}^2/\psi\le c_4$ on 
$[0,\infty)$ for some constant $c_4>0$. Let $R\ge 1$ and 
$$
\phi (x)=\psi\left (\frac{r(x)}{R}\right).
$$
Then
\begin{equation}
\frac{|\nabla\phi|^2}{\phi}=\frac{|\nabla^tr(x)|^2}{R^2}
\cdot\frac{\psi'{}^2}{\psi}\le\frac{c_3c_4}{R^2}.
\end{equation}
By an argument similar to \cite{C} and \cite{LY} we may assume without
loss of generality that $\phi (x)$ and $r(x)$ are both smooth on $M$. 
By (2.10), (2.11), the Laplacian comparison theorem \cite{SY} and an 
argument similar to the proof of Lemma 1.3 of \cite{Hs1} there exists 
a constant $c_5>0$ such that
\begin{align}
\Delta^tr(x)\le&c_5\quad\forall x\in M, r(x)\ge R, 0\le t\le T
\nonumber\\
\Rightarrow\qquad\Delta^t\phi=&\psi''\frac{|\nabla^tr(x)|^2}{R^2}+\psi'
\frac{\Delta^tr(x)}{R}
\ge-\frac{c_3c_4}{R^2}-\frac{c_4c_5}{R}
\ge -C_6(R^{-1}+R^{-2})
\end{align}
in $M\times [0,T]$ for some constant $C_6>0$. Suppose the function
$\phi F$ attains its maximum at $(x_0,t_0)$. If $(\phi F)(x_0,t_0)\le 0$
for any $R\ge 1$, then (2.1) holds in $M\times [0,T]$ and we are done.
Hence we may assume without loss of generality that there exists $R_1\ge 1$
such that for any $R\ge R_1$, $(\phi F)(x_0,t_0)>0$.
Then at $(x_0,t_0)$, $F_t\ge 0$,
\begin{equation}
\nabla (\phi F)=0\quad\Rightarrow\quad\phi\nabla F=-F\nabla\phi
\end{equation}
and by (2.2), (2.9), (2.12), (2.13) and (2.14),
\begin{align}
0\ge&\Delta (\phi F)\nonumber\\
\ge&\phi\Delta F+2\nabla\phi\cdot\nabla F+F\Delta\phi\nonumber\\
\ge&t_0\phi\frac{2(1-\delta)}{n}(\Delta f)^2-C_3t_0\phi
|\nabla f|^2-C_4t_0\phi-\frac{F}{t_0}\phi+a\phi F+\phi F_t
+2F\nabla f\cdot\nabla\phi\nonumber\\
&\qquad -2F\frac{|\nabla\phi|^2}{\phi}-C_6F(R^{-1}+R^{-2})\nonumber\\
\ge&t_0\phi\frac{2(1-\delta)}{n}(|\nabla f|^2-(f_t+af+q))^2-C_3t_0\phi
|\nabla f|^2-C_4t_0\phi-\frac{F}{t_0}\phi-a_-\phi F\nonumber\\
&\qquad -C_5F|\nabla f|\phi^{\frac{1}{2}}R^{-1}-C_6F(R^{-1}+R^{-2})
\end{align}
for some constant $C_5>0$. Let
$$
\sigma=\frac{\alpha^2-1}{\alpha^2},\quad 
\delta_1=\frac{2(\alpha -1)}{\alpha^2}\quad
\mbox{ and }\tau=\frac{1}{2}(\sigma -\delta_1).
$$
Then by (2.15) and an argument similar to the proof on P.13--14 of \cite{CTY},
\begin{align}
0\ge&-\phi F[C_7t_0(1+\delta_1^{-1})(R^{-1}+R^{-2})+1+t_0a_-]-C_8(1+\tau^{-1})t_0^2
\nonumber\\
&\qquad
+\frac{2(1-\delta)}{n}(1-\sigma)(t_0\phi)^2(|\nabla f|^2-\alpha (f_t+af+q))^2
\nonumber\\
\ge&-\phi F[C_7t_0(1+\delta_1^{-1})(R^{-1}+R^{-2})+1+t_0a_-]
-C_8(1+\tau^{-1})t_0^2
+\frac{2(1-\delta)}{n}\alpha^{-2}(\phi F)^2
\end{align}
for some constant $C_7>0$ and $C_8>0$. Let $y=\phi F$ and $b=
C_7t_0(1+\delta_1^{-1})(R^{-1}+R^{-2})+1+t_0a_-$. Then by (2.16),
\begin{align}
&y^2-\frac{n\alpha^2b}{2(1-\delta)}y\le\frac{C_8n}{2(1-\delta)}
\alpha^2(1+\tau^{-1})t_0^2\nonumber\\
\Rightarrow\quad&\left(y-\frac{n\alpha^2b}{4(1-\delta)}\right)^2
\le\left(\frac{n\alpha^2b}{4(1-\delta)}\right)^2
+\frac{C_8n}{2(1-\delta)}\alpha^2(1+\tau^{-1})t_0^2\nonumber\\
\Rightarrow\quad&y\le\frac{n\alpha^2}{2(1-\delta)}
[C_7T(1+\delta_1^{-1})(R^{-1}+R^{-2})+1+Ta_-]
+C_9T.
\end{align}
where $C_9=\alpha\left(\frac{C_8n}{2(1-\delta)}\right)^{\frac{1}{2}}
(1+\tau^{-1})^{\frac{1}{2}}$. Then for any $x\in B_R(x_0)$,
$$
F(x,T)\le\frac{n\alpha^2}{2(1-\delta)}[C_7T(1+\delta_1^{-1})(R^{-1}+R^{-2})
+1+Ta_-]+C_9T
$$
Replacing $T$ by $t$ and letting $t\in (0,T]$ we have for any 
$x\in B_R(x_0)$,
\begin{align}
&|\nabla f|^2-\alpha f_t-\alpha af-\alpha q
\le\frac{n\alpha^2}{2(1-\delta)}[C_7(1+\delta_1^{-1})(R^{-1}+R^{-2})+
\frac{1}{t}+a_-]+C_9\nonumber\\
\Rightarrow\quad&|\nabla f|^2-\alpha f_t-\alpha af-\alpha q
\le\frac{n\alpha^2}{2(1-\delta)t}(1+ta_-)
+C_9\quad\mbox{ as }R\to\infty\nonumber.
\end{align}
and the theorem follows.
\end{proof}

\begin{cor}
Let $(M, g(t))$, $0\le t\le T$, be a complete noncompact manifold with metric 
$g$ evolving by the Ricci flow (1.1) and 
$$
|\nabla^iRm|\le k_0\quad\text{ on }M\times [0,T]\quad\forall i=0,1
$$
for some constant $k_0>0$. Let $u$ be a positive solution of (1.5). Then 
for any 
$\alpha>1$ and $0<\delta<1$ there exists a constant $C_1>0$ depending on 
$k_0$, $\alpha$, $\delta$, $a_+$, the space-time uniform bound of 
$|q|$, $|\nabla^tq|$, $|\Delta^tq|$ such that
\begin{equation*}
\frac{|\nabla u|^2}{u^2}-\alpha\frac{u_t}{u}-\alpha a\log u-\alpha q
\le\frac{n\alpha^2}{2(1-\delta)t}(1+ta_-)+C_1\quad\text{ in }M\times (0,T].
\end{equation*}
\end{cor}

\begin{cor}
Let $(M, g)$ be a complete noncompact manifold with a fixed metric 
$g$ satisfying
$$
|Rm|\le k_0\quad\text{ in }M
$$
for some constant $k_0>0$. Let $u$ be a positive solution of 
$$
\Delta u-au\log u-qu=0\quad\text{ in }M
$$
where $q$ is a smooth function on $M$ with $|q|$, $|\nabla q|$
and $|\Delta q|$ uniformly bounded on $M$. Then there exists a  constant
$C_2>0$ such that 
\begin{equation*}
u\le C_2\quad\mbox{ in }M\quad\mbox{ if }a<0
\end{equation*}
and there exists a  constant $C_3>0$ such that 
\begin{equation*}
u\ge C_3\quad\mbox{ in }M\quad\mbox{ if }a>0.
\end{equation*}
\end{cor}
\begin{proof}
By Theorem 2.1 with $h_{ij}=0$, $\alpha=2$, $\delta=1/2$, there exists a 
constant $C_1>0$ such that (2.1) holds for all $t>0$. Letting $t\to\infty$
in (2.1),
\begin{align*}
&\frac{|\nabla u|^2}{u^2}-2a\log u-2q\le 4na_-+C_1\nonumber\\
\Rightarrow\quad &-2a\log u-2q\le 4na_-+C_1\nonumber\\
\Rightarrow\quad &\left\{\begin{aligned}
&u\le e^{(4na_-+C_1+2\|q\|_{\infty})/2|a|}\quad\mbox{ in }M
\quad\mbox{ if }a<0\\
&u\ge e^{-(C_1+2\|q\|_{\infty})/2a}\qquad\quad\mbox{ in }M
\quad\mbox{ if }a>0\end{aligned}\right.
\end{align*}
and the corollary follows.
\end{proof}

\end{document}